\documentclass[12pt]{article}

\usepackage{amssymb}
\usepackage{epsfig}
\usepackage{amsmath}
\usepackage{amsthm}
\usepackage{subfigure}
\usepackage{url}
\usepackage{color}
\usepackage{graphicx}
\usepackage[margin=1in]{geometry}
\DeclareMathOperator{\transverse}{\cap\kern-7.75pt\top}

\newcommand{\codim}{\operatorname{codim}}

\newcommand{\R}{\Bbb{R}}

\newtheorem{exm}{Example}

\usepackage[bookmarks=true, % true means bookmarks
            bookmarksnumbered=false, % left window are numbered
            bookmarksopen=false, % true means only level 1 are displayed. 
            colorlinks=true, 
            linkcolor=red]{hyperref}

\begin{document}

\title{{\bf Some Minimal Shape Decompositions Are Nice}} 
%\title{{\bf Simplicity is Conserved by Some Shape Decompositions}}
\author{Kevin R. Vixie} 
\date{Washington State University}

\maketitle

\subsection*{In a Nutshell}
\label{synopsis}

In some sense, the world is composed of \emph{shapes} and
\emph{words}, of continuous things and discrete things.  The
recognition and study of continuous objects in the form of shapes
occupies a significant part of the effort of unraveling many
geometric questions. Shapes can be represented with great generality
by objects called \emph{currents}.  While the enormous variety and
representational power of currents is useful for representing a huge
variety of phenomena, it also leads to the problem that knowing
something is a respectable current tells you little about how nice or
\emph{regular} it is. In these brief notes I give an intuitive
explanation of a result that says that an important class of minimal
shape decompositions will be nice if the input shape (current) is
nice. \emph{\color{blue}These notes are an exposition of the paper by
  Ibrahim, Krishnamoorthy and Vixie~\cite{ibrahim-2014-flat-arXiv}
  which can be found on the arXiv:
  \href{http://arxiv.org/abs/1411.0882}{http://arxiv.org/abs/1411.0882}
  and any reference to these notes, should include a reference to that
  paper as well.}

\subsection*{Shapes and Currents}
\label{sec:shapes-and-currents}

% In some sense, the world is composed of \emph{shapes} and
% \emph{words}, of continuous analytic/geometric things and discrete
% algebraic/categorical/linguistic things.  The recognition and study of
% shapes forms a significant part of the effort of unravelling many
% geometric questions.

We begin by representing shapes as k-\emph{currents} which, for this
section, should be thought of as oriented k-dimensional subsets of
$\R^n$ ($0\leq k \leq n$) that are locally flat almost everywhere\footnote{By
  \emph{locally flat almost everywhere}, we mean that when we zoom in to almost any
  point in the current, it looks (measure theoretically) like a
  k-dimensional plane. Such sets are called \emph{rectifiable sets}.}.

It turns out that currents\footnote{To be more precise, k-currents are
  defined to be elements of the dual space to the space of smooth
  k-forms in $\R^n$. But for the purposes of these notes, thinking of
%%  TLS edit
%%  the original line is commented out
% them as unions of of pieces of oriented k-manifolds will get you a
  them as unions of pieces of oriented k-manifolds will get you a
  long way towards understanding currents.} are a deep generalization of
sets like 2-dimensional surfaces in $\R^3$ or 1-dimensional curves in
$\R^2$ or $\R^3$. This is both very helpful, because we can represent
a large class of things very naturally, and challenging, because
knowing something is a respectable current no longer implies it is in
some small class of relatively tame objects.

As a result of this diversity, we are naturally led to classify
currents into different subspaces -- e.g. integral,
normal, or flat currents -- and when we find that some problem has a
current as a solution, we then occupy ourselves with finding out which
of these classes that current falls into. This in turn tells us how
tame (how \emph{regular}) it is.  The nicest class is the class of
integral currents. Proving that a current is in a smaller, more
regular set of currents can sometimes be very challenging.

\subsection*{Shapes and Currents: Examples}
\label{sec:sc-examples}

\begin{exm}[Orientation]
  Suppose you have a simple closed curve in $\R^2$ -- it does not
  cross itself and is really just a circle that has (possibly) been
%%  TLS edit
%%  the original line is commented out
% contorted stretched around. Then an orientation is a choice of
  contorted or stretched around. Then an orientation is a choice of
  direction along the curve: you can choose a unit tangent vector
  everywhere such that the directions are all consistent. Notice that
%%  TLS edit
%%  the original line is commented out
% there will be two possible choices, clockwise and counter-clockwise.
  there will be two possible choices, clockwise and counterclockwise.
  For a closed 2-dimensional surface in $\R^3$ -- say a sphere is
  $\R^3$ -- we can orient the surface by choosing little tangent
  planes with normal vectors pointing either inside or outside of the
  surface. Even if the surface is not closed, we can attempt to assign
%%  TLS edit
%%  the original line is commented out
% a the little tangent patches and normal vector to every point in a
  a little tangent patch and normal vector to every point in a
  way that is continuous. As long as all these choices are consistent
%%  TLS edit
%%  the original line is commented out
% (the normal vectors change continuously) then we cay the surface is
  (the normal vectors change continuously) then we say the surface is
  oriented by the little 2-planes which we call \emph{2-vectors}.
\begin{figure}[htp!]
\begin{center}
\subfigure[Oriented Curve]{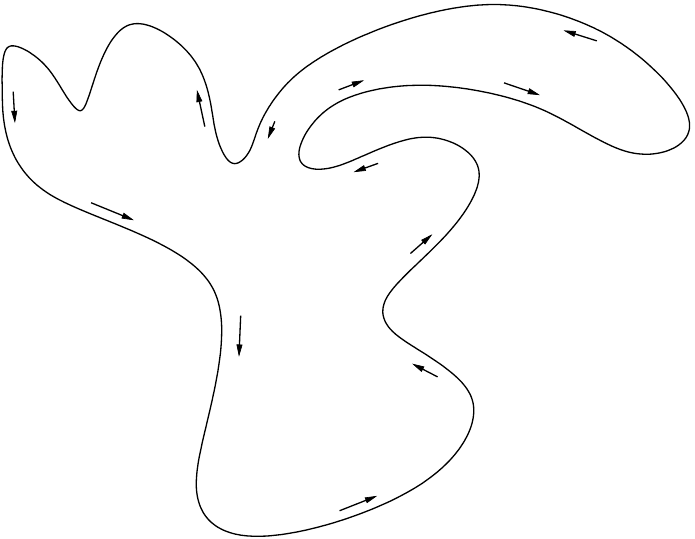}
\subfigure[Oriented Surface]{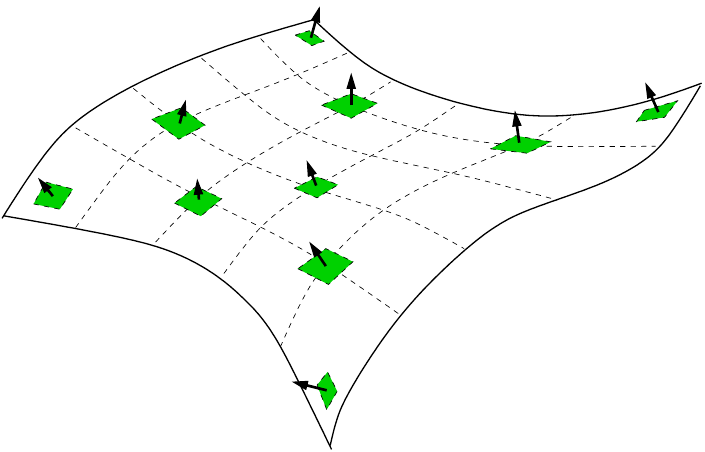}
\end{center}
\end{figure}
\end{exm}
\begin{exm}[Boundaries]
%%  TLS edit
%%  the original line is commented out
% We will denote the boundary of $T$ by $\partial T$. If $T$ is a the
  We will denote the boundary of $T$ by $\partial T$. If $T$ is the
  2-current in $\R^2$ corresponding to the unit disk with
%%  TLS edit
%%  the original line is commented out
% counterclockwise (right-handed) orientation, the $\partial T$ is the
  counterclockwise (right-handed) orientation, $\partial T$ is the
%%  TLS edit
%%  the original line is commented out
% unit circle oriented in a counter-clockwise direction.
  unit circle oriented in a counterclockwise direction.
\begin{center}
\begin{picture}(0,0)%
\includegraphics{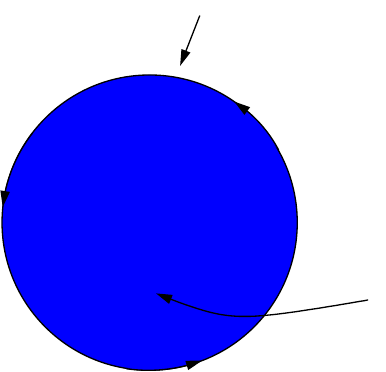}%
\end{picture}%
\setlength{\unitlength}{1657sp}%
\begingroup\makeatletter\ifx\SetFigFont\undefined%
\gdef\SetFigFont#1#2#3#4#5{%
  \reset@font\fontsize{#1}{#2pt}%
  \fontfamily{#3}\fontseries{#4}\fontshape{#5}%
  \selectfont}%
\fi\endgroup%
\begin{picture}(4375,4211)(900,-4146)
\put(3360,-202){\makebox(0,0)[lb]{\smash{{\SetFigFont{8}{9.6}{\rmdefault}{\mddefault}{\updefault}{\color[rgb]{0,0,0}$\partial T$}%
}}}}
\put(5260,-3377){\makebox(0,0)[lb]{\smash{{\SetFigFont{8}{9.6}{\rmdefault}{\mddefault}{\updefault}{\color[rgb]{0,0,0}$T$}%
}}}}
\end{picture}%

\end{center}
\end{exm}

\begin{exm}[Mass of a current] $M(T)$ denotes the \emph{mass} of the
  k-current, which for the purposes of these notes should be thought
%%  TLS edit
%%  the original line is commented out
% of as the either the k-dimensional volume of the current or as the
  of as either the k-dimensional volume of the current or as the
  multiplicity weighted k-dimensional volume of the current.  For
%%  TLS edit
%%  the original line is commented out
% example, If $T$ is a the 2-current in $\R^2$ corresponding to the
  example, If $T$ is the 2-current in $\R^2$ corresponding to the
%%  TLS edit
%%  the original line is commented out
% unit disk then $M(T) = \pi$ and the mass of $M(\partial T) =
  unit disk then $M(T) = \pi$ and the mass $M(\partial T) =
  2\pi$. In the important case in which we have $T_\Omega$, the
  current corresponding to an oriented surface $\Omega$, we get that
  $M(T_\Omega) = \int_\Omega 1 dx$ as is illustrated here.
\begin{center}
\begin{picture}(0,0)%
\includegraphics{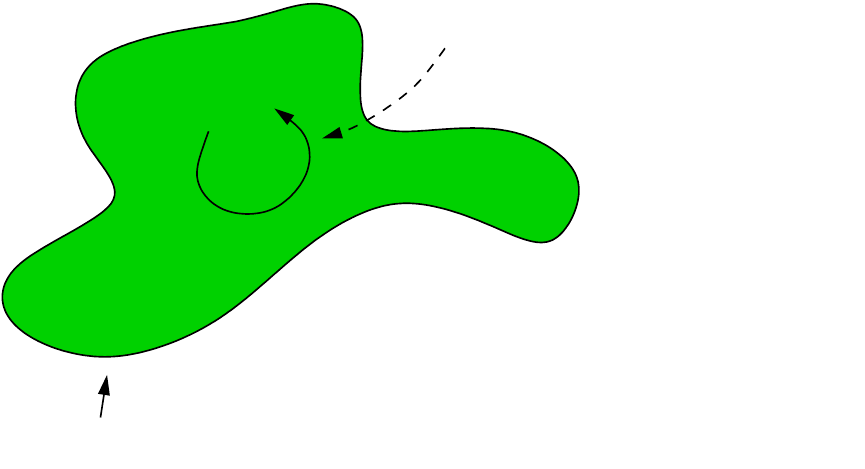}%
\end{picture}%
\setlength{\unitlength}{2072sp}%
\begingroup\makeatletter\ifx\SetFigFont\undefined%
\gdef\SetFigFont#1#2#3#4#5{%
  \reset@font\fontsize{#1}{#2pt}%
  \fontfamily{#3}\fontseries{#4}\fontshape{#5}%
  \selectfont}%
\fi\endgroup%
\begin{picture}(7889,4274)(1153,-5342)
\put(1667,-3928){\makebox(0,0)[lb]{\smash{{\SetFigFont{10}{12.0}{\rmdefault}{\mddefault}{\updefault}{\color[rgb]{0,0,0}$\Omega$}%
}}}}
\put(1457,-5233){\makebox(0,0)[lb]{\smash{{\SetFigFont{10}{12.0}{\rmdefault}{\mddefault}{\updefault}{\color[rgb]{0,0,0}$T_\Omega$}%
}}}}
\put(4667,-4438){\makebox(0,0)[lb]{\smash{{\SetFigFont{10}{12.0}{\rmdefault}{\mddefault}{\updefault}{\color[rgb]{0,0,0}$M(T_\Omega) = \int_\Omega 1 dx$}%
}}}}
\put(5012,-1318){\makebox(0,0)[lb]{\smash{{\SetFigFont{10}{12.0}{\rmdefault}{\mddefault}{\updefault}{\color[rgb]{0,0,0}Counterclockwise orientation}%
}}}}
\end{picture}%

\end{center}
\end{exm}

\begin{exm}[Adding and Subtracting Currents]
  A key feature of currents is the fact that \emph{it makes sense to
    add and subtract currents}: Suppose that $T_1$ is the unit disk
%%  TLS edit
%%  the original line is commented out
% $T_1$ oriented clockwise and $T_2$ is the intersection
  $T_1$ oriented counterclockwise and $T_2$ is the intersection
%%  TLS edit
%%  the original line is commented out
  of the counterclockwise oriented unit disk intersected with the upper half
  plane. Then $T_3 = T_1 + T_2$ will be the 2-dimensional current
  corresponding to the unit disk with multiplicity 2 in the upper half
  disk, and multiplicity 1 in the lower half disk. It's boundary will
  be the union of the circle of unit radius whose upper half has
  multiplicity 2 and lower half has multiplicity 1 and the diameter
  segment through the origin with multiplicity 1.
\begin{center}
\begin{picture}(0,0)%
\includegraphics{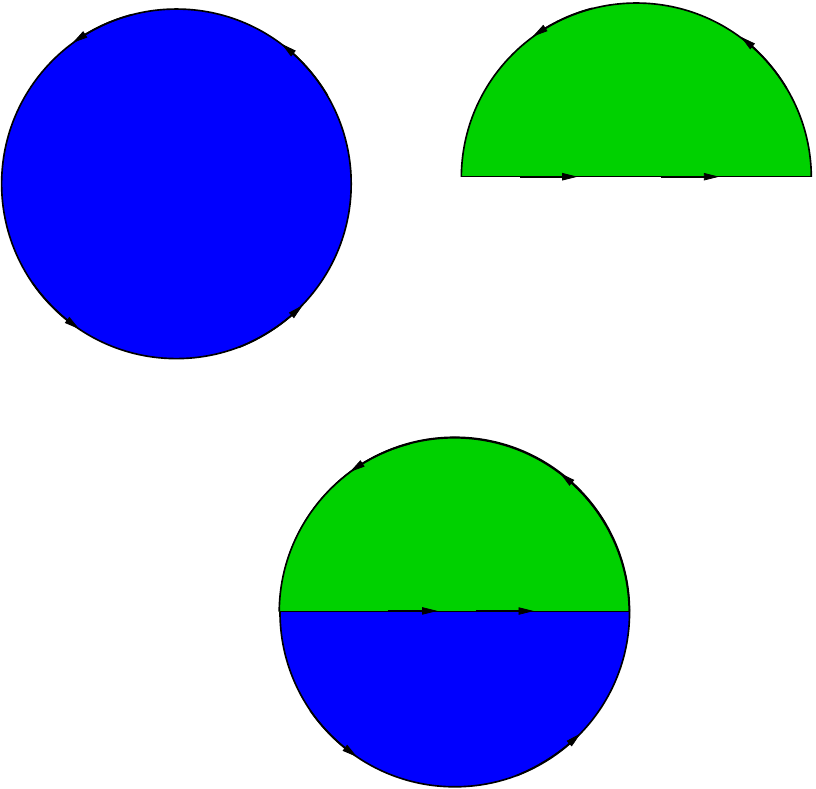}%
\end{picture}%
\setlength{\unitlength}{2072sp}%
\begingroup\makeatletter\ifx\SetFigFont\undefined%
\gdef\SetFigFont#1#2#3#4#5{%
  \reset@font\fontsize{#1}{#2pt}%
  \fontfamily{#3}\fontseries{#4}\fontshape{#5}%
  \selectfont}%
\fi\endgroup%
\begin{picture}(7436,7197)(547,-6956)
\end{picture}%

\end{center}
Subtracting $T_2$ from $T_1$, we simply get the lower half unit disk or
lower semi-disk, again oriented counterclockwise.
\end{exm}
\begin{exm}[Mass, again]
%%  TLS edit
%%  the original line is commented out
% If $T_1$ is the clockwise oriented unit disk at the origin, and
  If $T_1$ is the counterclockwise oriented unit disk at the origin, and
%%  TLS edit
%%  the original line is commented out
% $T_2$ is the clockwise oriented unit disk centered at $(0,3)$, then
  $T_2$ is the counterclockwise oriented unit disk centered at $(0,3)$, then
  $T = 2T_1 + 3T_2$, is the sum of the oriented unit disk
  centered at the origin with multiplicity $2$ and oriented the unit
  disk centered at $(0,3)$ with multiplicity $3$. The mass of $T$ is
%%  TLS edit
%%  the original line is commented out
% $5\pi$ and the $M(\partial T)$ is $10\pi$. In general, if the
  $5\pi$ and $M(\partial T)$ is $10\pi$. In general, if the
  multiplicity of the current is $N(x)$ then $M(T_\Omega) =
  \int_\Omega N(x) dx$.
\begin{center}
\begin{picture}(0,0)%
\includegraphics{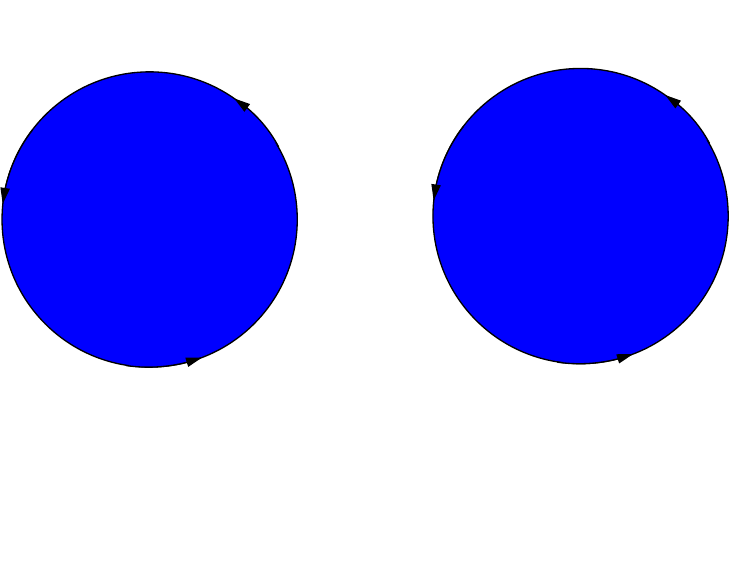}%
\end{picture}%
\setlength{\unitlength}{1657sp}%
\begingroup\makeatletter\ifx\SetFigFont\undefined%
\gdef\SetFigFont#1#2#3#4#5{%
  \reset@font\fontsize{#1}{#2pt}%
  \fontfamily{#3}\fontseries{#4}\fontshape{#5}%
  \selectfont}%
\fi\endgroup%
\begin{picture}(8339,6690)(900,-6641)
\put(1935,-5212){\makebox(0,0)[lb]{\smash{{\SetFigFont{12}{14.4}{\rmdefault}{\mddefault}{\updefault}{\color[rgb]{0,0,0}$M(T) = M(2T_1 + 3T_2)$}%
}}}}
\put(3185,-6487){\makebox(0,0)[lb]{\smash{{\SetFigFont{12}{14.4}{\rmdefault}{\mddefault}{\updefault}{\color[rgb]{0,0,0}$= 5\pi$}%
}}}}
\put(3185,-5837){\makebox(0,0)[lb]{\smash{{\SetFigFont{12}{14.4}{\rmdefault}{\mddefault}{\updefault}{\color[rgb]{0,0,0}$= 2M(T_1) + 3M(T_2)$}%
}}}}
\put(3035,-302){\makebox(0,0)[lb]{\smash{{\SetFigFont{12}{14.4}{\rmdefault}{\mddefault}{\updefault}{\color[rgb]{0,0,0}$T = 2*T_1 + 3*T_2$}%
}}}}
\put(7185,-2552){\makebox(0,0)[lb]{\smash{{\SetFigFont{12}{14.4}{\rmdefault}{\mddefault}{\updefault}{\color[rgb]{0,0,0}$3T_2$}%
}}}}
\put(2160,-2577){\makebox(0,0)[lb]{\smash{{\SetFigFont{12}{14.4}{\rmdefault}{\mddefault}{\updefault}{\color[rgb]{0,0,0}$2T_1$}%
}}}}
\end{picture}%

\end{center}
\end{exm}

\subsection*{Distances Between Shapes: The Multiscale Flat Norm}
\label{sec:msfn}

A first requirement for studying collections of objects is a distance
between those objects. When those objects are currents there is a very
natural distance called the \emph{flat norm}, $\Bbb{F}_{\lambda}(T)$. This
%%  TLS edit
%%  the original line is commented out
 % distance begins be decomposing the input k-current $T$ into two pieces
distance begins by decomposing the input k-current $T$ into two pieces
$T = (T -\partial S) + \partial S$, where $S$ is a k+1-current. Now we
%%  TLS edit
%%  the original line is commented out
 % measure the cost of the decomposition is $M(T -\partial S) + \lambda M(S)$ .
measure the cost of the decomposition as $M(T -\partial S) + \lambda M(S)$ .
Finally, we minimize this cost over all possible decompositions
(i.e. over all possible k+1 currents $S$):
\[ \Bbb{F}_{\lambda}(T) = \min_S\left( M(T -\partial S) + \lambda M(S)\right).\] 
\begin{exm}[The multiscale flatnorm] In this example, a 1-current $T$
  is decomposed into two pieces, $T-\partial S$ and $\partial S$, in a way that minimizes
  \[M(T-\partial S) + \lambda M(S). \] It turns out that in this case
  it can be shown that the optimal decomposition fills in the corners
  of $T$ with $S$ such that the resulting $T-\partial S$ has curvature
  bounded by $\frac{1}{\lambda}$; that is, we short-cut the corners
  with arcs of circles of radius $\frac{1}{\lambda}$.
\begin{center}
\begin{picture}(0,0)%
\includegraphics{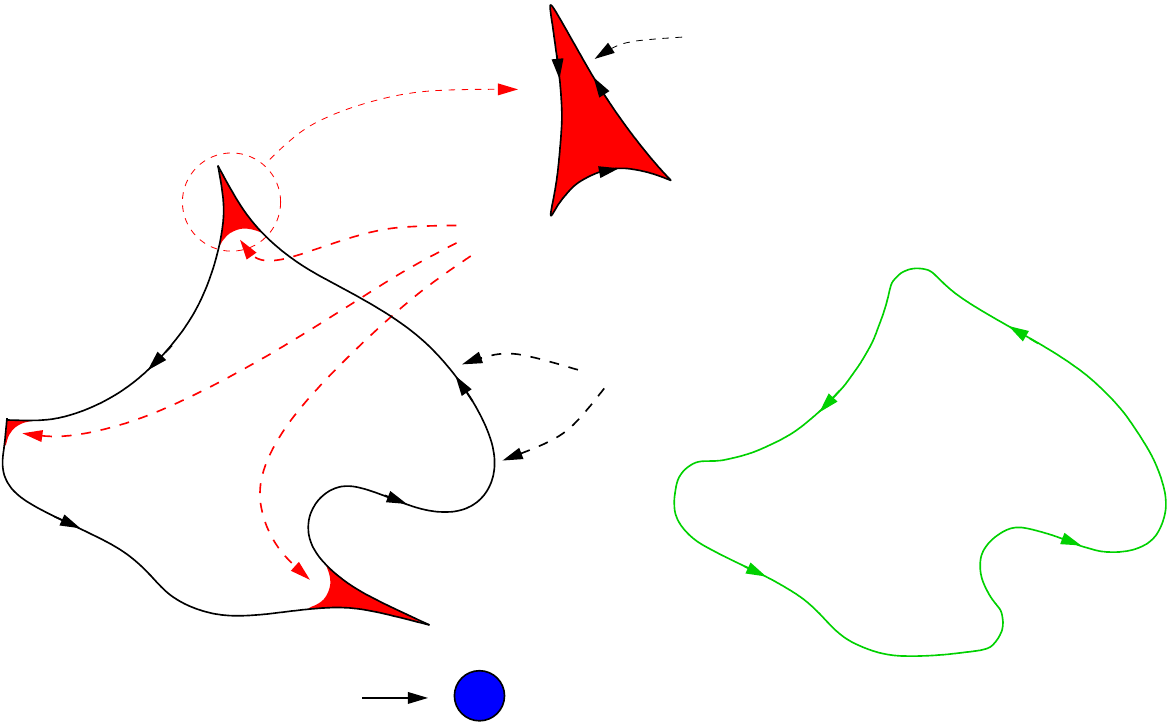}%
\end{picture}%
\setlength{\unitlength}{2072sp}%
\begingroup\makeatletter\ifx\SetFigFont\undefined%
\gdef\SetFigFont#1#2#3#4#5{%
  \reset@font\fontsize{#1}{#2pt}%
  \fontfamily{#3}\fontseries{#4}\fontshape{#5}%
  \selectfont}%
\fi\endgroup%
\begin{picture}(10680,6581)(1663,-5598)
\put(3795,-5459){\makebox(0,0)[lb]{\smash{{\SetFigFont{10}{12.0}{\rmdefault}{\mddefault}{\updefault}{\color[rgb]{0,0,0}$\frac{1}{\lambda}$ ball}%
}}}}
\put(5955,-1234){\makebox(0,0)[lb]{\smash{{\SetFigFont{10}{12.0}{\rmdefault}{\mddefault}{\updefault}{\color[rgb]{1,0,0}$S$}%
}}}}
\put(7005,-2444){\makebox(0,0)[lb]{\smash{{\SetFigFont{10}{12.0}{\rmdefault}{\mddefault}{\updefault}{\color[rgb]{0,0,0}$T$}%
}}}}
\put(8920,-3589){\makebox(0,0)[lb]{\smash{{\SetFigFont{10}{12.0}{\rmdefault}{\mddefault}{\updefault}{\color[rgb]{0,.82,0}$T-\partial S$}%
}}}}
\put(7967,588){\makebox(0,0)[lb]{\smash{{\SetFigFont{10}{12.0}{\rmdefault}{\mddefault}{\updefault}{\color[rgb]{0,0,0}$\partial S$}%
}}}}
\end{picture}%

\end{center}
\end{exm}
While it might seem difficult to compute the flat norm, the
realization that the flat norm is connected to the $L^1$TV functional
from image analysis\cite{morgan-2007-1} permits efficient calculation
for the case when $T$ is an (n-1)-dimensional boundary of an
n-dimensional set in $\R^n$. (In the figure above, $T$ is such an
%%  TLS edit
%%  the original line is commented out
 % object: it is the 1 dimensional boundary of a 2 dimensional subset of
object: it is the 1-dimensional boundary of a 2-dimensional subset of
the plane.) In this ``co-dimension 1 boundary'' case, we can use a
variety of algorithms to find the minimal decompositions. One efficient
method is the graph-cut method which has been used to to calculate
good approximations to the minimal decompositions for shapes in
$\R^2$\cite{vixie-2010-2}.
\begin{exm}[Differences between Shapes] Now we are in a position to
  consider distances between shapes. In the figure, we notice that
  while a direct subtraction of $T_1$ from $T_2$ produces something
  whose mass is not small due to the fact that they do not coincide
  and cancel, the flat norm still sees the two currents as very
%%  TLS edit
%%  the original line is commented out
% close and assigns a small distance to the difference between $T_1$
  close and assigns a small distance between $T_1$
  and $T_2$.
\begin{center}
\begin{picture}(0,0)%
\includegraphics{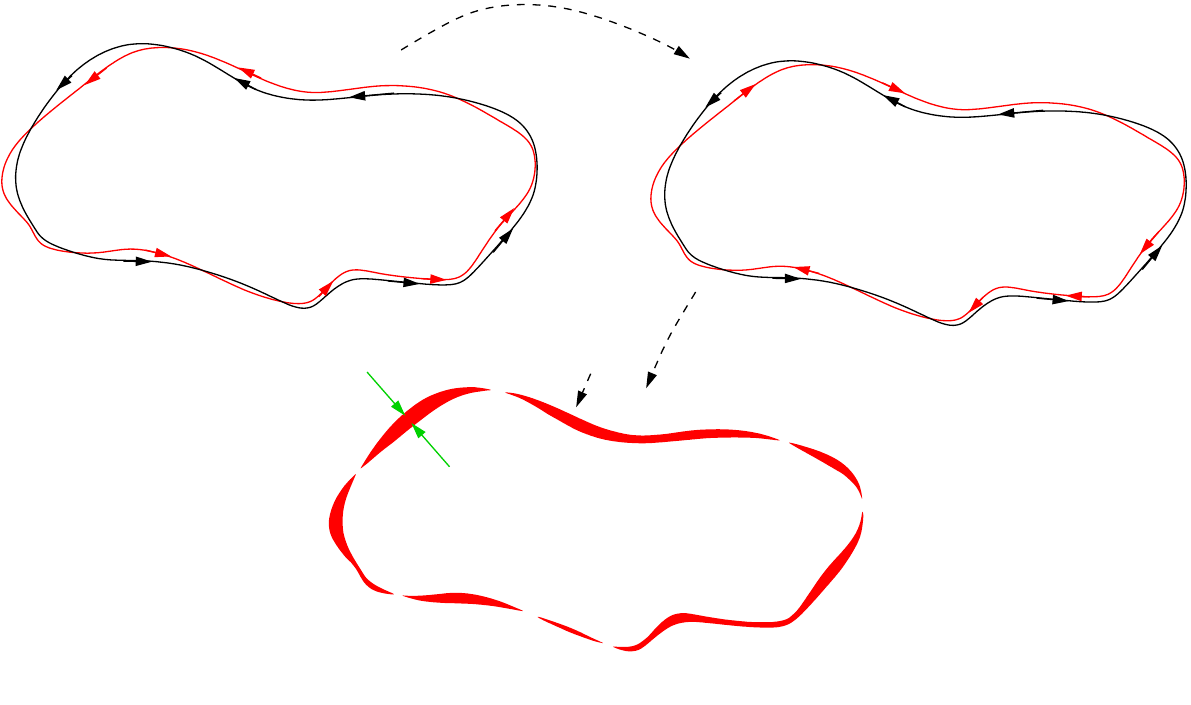}%
\end{picture}%
\setlength{\unitlength}{1657sp}%
\begingroup\makeatletter\ifx\SetFigFont\undefined%
\gdef\SetFigFont#1#2#3#4#5{%
  \reset@font\fontsize{#1}{#2pt}%
  \fontfamily{#3}\fontseries{#4}\fontshape{#5}%
  \selectfont}%
\fi\endgroup%
\begin{picture}(13580,8025)(1301,-8748)
\put(4106,-1335){\makebox(0,0)[lb]{\smash{{\SetFigFont{8}{9.6}{\rmdefault}{\mddefault}{\updefault}{\color[rgb]{1,0,0}$T_2$}%
}}}}
\put(3489,-1952){\makebox(0,0)[lb]{\smash{{\SetFigFont{8}{9.6}{\rmdefault}{\mddefault}{\updefault}{\color[rgb]{0,0,0}$T_1$}%
}}}}
\put(10472,-2869){\makebox(0,0)[lb]{\smash{{\SetFigFont{8}{9.6}{\rmdefault}{\mddefault}{\updefault}{\color[rgb]{0,0,0}$T_1-T_2$}%
}}}}
\put(6452,-6747){\makebox(0,0)[lb]{\smash{{\SetFigFont{8}{9.6}{\rmdefault}{\mddefault}{\updefault}{\color[rgb]{1,0,0}$\partial S = T_1 - T_2$}%
}}}}
\put(4490,-4821){\makebox(0,0)[lb]{\smash{{\SetFigFont{8}{9.6}{\rmdefault}{\mddefault}{\updefault}{\color[rgb]{0,.82,0}thickness $< < 1$}%
}}}}
\put(10745,-4985){\makebox(0,0)[lb]{\smash{{\SetFigFont{8}{9.6}{\rmdefault}{\mddefault}{\updefault}{\color[rgb]{0,0,0}$M(T_1 - T_2) = M(T_1) + M(T_2)$ is not small!}%
}}}}
\put(4878,-8639){\makebox(0,0)[lb]{\smash{{\SetFigFont{8}{9.6}{\rmdefault}{\mddefault}{\updefault}{\color[rgb]{1,0,0}$\Bbb{F}_\lambda(T_1-T_2) = M(T_1-T_2-\partial S) + \lambda M(S) = \lambda M(S)$ is small!}%
}}}}
\end{picture}%

\end{center}
\end{exm}
Now we are ready for the main result in the paper we are aiming to explain.
\subsection*{When does Integral T $\rightarrow$ Integral S?}
\label{sec:intT-intS}

% Now we come to the central question addressed by the paper we are
% discussing: how nice are the flat norm minimizers. In our case, we
% define nice to be any current that is an \emph{integral
%   current}. These correspond to unions of pieces of nice k-dimensional
% surfaces in $\R^n$.

% {\color{blue}If $T$ is integral, will the optimal
%   $S$ (and therefore the optimal \{$T -\partial S$, $\partial S$\}
%   decomposition) also be integral?}.

% We start with a definition.

% \begin{deff}[Integral Current]

% In our
% case, we define nice to be any current that is an \emph{integral
%   current}. These correspond to unions of pieces of nice k-dimensional
% surfaces in $\R^n$. {\color{blue}If $T$ is integral, will the optimal
%   $S$ (and therefore the optimal \{$T -\partial S$, $\partial S$\}
%   decomposition) also be integral?}.
  
% \end{deff}

Using currents to represent shapes, and the flat norm to measure size and distances in the 
space of currents leads us to the question:

\bigskip

{\color{blue}Is the $S$ minimizing $M(T-\partial S) + \lambda M(S)$ nice when the input T is nice? }

\bigskip

\noindent In our case, we consider a current to be nice if it is an
\emph{integral current}. These correspond to unions of pieces of
smooth\footnote{Actually, by smooth we mean that it is $C^1$, that is
  the surface is a submanifold of $\R^n$ whose tangents vary
  continuously.} k-dimensional surfaces in $\R^n$, such that the
%%  TLS edit
%%  the original line is commented out
 % k-dimensional volume of the union and the k-1 dimensional volume of
k-dimensional volume of the union and the (k-1)-dimensional volume of
its boundary are both finite. So the question can be restated:

\bigskip

{\color{blue}Is the $S$ minimizing $M(T-\partial S) + \lambda M(S)$ integral when T is integral? }

\bigskip

\noindent To get to the answer, we introduce a bit more terminology:
We say a set $E\subset\R^n$ has co-dimension $k$ if $\dim(E) = n-k.$
\begin{exm}[co-dimension]
   A 1-dimensional curve in $\R^2$ has co-dimension 1, whereas a 1-dimensional
curve in $\R^3$ has co-dimension 2.
\end{exm}

Previous to this work what was known could be summed up in two statements:

\begin{enumerate}
\item Together~\cite{chan-2005-1,morgan-2007-1} show that when $T$ is a co-dimension 1 boundary in
$\R^n$, there is a minimizing integral S\footnote{We
do not have uniqueness -- there are times in which there are also
non-integral minimizers along with the integral minimizer.}. 
\item In~\cite{ibrahim-2013-simplicial} it is shown that co-dimension 1 T's
have minimizing integral S's in the case that the T and S live on a
\emph{simplicial complex}. Note that T need not be a boundary, though it must be
a discrete, simplicial current!
\end{enumerate}

\subsection*{The Punchline}
\label{sec:punch}

The main result of the paper by Ibrahim, Krishnamoorthy and
Vixie\cite{ibrahim-2014-flat-arXiv} uses the results for discrete,
simplicial currents to prove that the flat norm decompositions of
arbitrary co-dimension 1 integral currents are also well behaved:

\begin{quote}{\color{blue} If T is a co-dimension 1 integral current,
    then the flat norm minimizing $S$ can also be taken to be
    integral. {\bf In the case of 1-currents in $\R^2$, the proof is
    complete.} In the case of co-dimension 1 currents in $\R^n$ for
%%  TLS edit
%%  the original line is commented out---is it better my way?
%   $n\geq 3$, the proof is complete provided the truth of a conjecture in
%   the paper concerning the existence of certain simplicial
%   refinements\footnote{We prove a case of the conjecture that gives
    $n\geq 3$, the proof is complete assuming a conjecture in
    the paper concerning the existence of certain simplicial
    refinements is true\footnote{We prove a case of the conjecture that gives
      us the result for 1-currents in $\R^2$. The full conjecture
      seems true, but may be difficult to prove.}.}
\end{quote}

In essence, the argument for the main result boils down to the
observation that if T is integral and there is not an integral
minimizer then this implies, after lots of fussing around and using
modified deformation theorems and special simplicial refinements and
adaptations, that we can produce a integral simplicial current that
also does not have an integral minimizer. And that is a contradiction,
according to the results in~\cite{ibrahim-2013-simplicial}.

\bigskip

\begin{quote}{\color{blue} What about T with $\codim(T) > 1$?}
\end{quote}

\bigskip

The paper also uses the isoperimetric inequality to prove that when
$T$ is a boundary and $\lambda$ is small enough, the flat norm
minimization becomes a minimal surface problem. This is used to
conclude that:

\bigskip

\begin{quote}{\color{blue} T being integral does not imply there is a minimizing S
  that is also integral when $\codim(T) \geq 3$}
\end{quote}

\bigskip

\noindent because of examples of 1-dimensional boundaries in $\R^4$ for which there
are spanning currents\footnote{$S$ such that $\partial S = T$} $S$ such
that $M(S) < M(E)$ for any integral $E$ with $\partial E =
T$~\cite{morgan-1984-area}.

%%  TLS edit
%%  the original line is commented out
 % The case of co-dimension 2 case is completely open, except for the
The case of co-dimension 2 is completely open, except for the
fact that we know we cannot use minimal surface problems to find
counterexamples~\cite{morgan-1984-area}.
\subsection*{Post Script}
\label{comment}
While it is true that integral currents are \emph{much} more regular
than arbitrary currents, they are \emph{not} almost smooth manifolds.
In order to dispel that notion, we give some examples of nice but
non-trivial integral 1-currents in $\R^2$.
\begin{figure}[htp!]
  \centering \subfigure[A perfectly nice integral 1-current. Since the
  boundary has to be finite, there must be at most a finite number of
  non-closed curves since each such curve contributes 2 to the mass of
  the boundary. ]{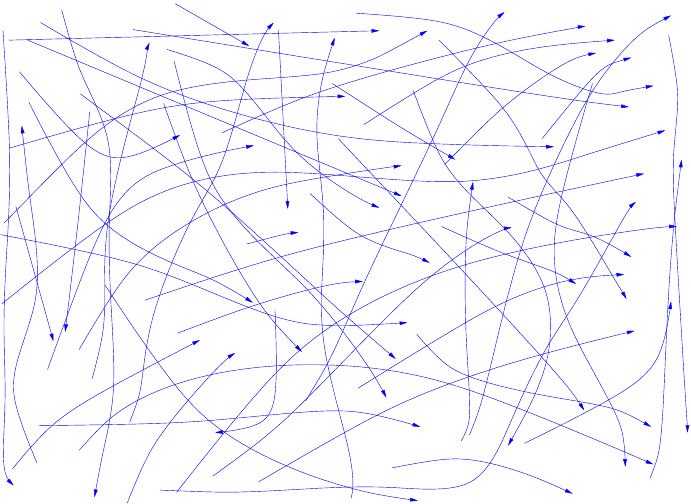} \subfigure[Another
  perfectly nice integral current. The circles are all oriented in the
  counterclockwise direction.]{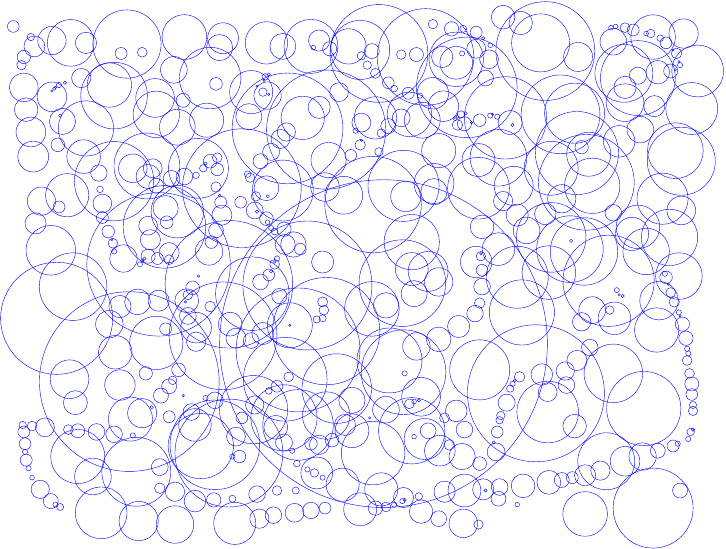}
\end{figure}
While it is true that there can only be a finite number of non-closed
curves in an integral 1-current, that doesn't mean such a 1-current is
simple. A current that is the union of 100,000 curves each with length
%%  TLS edit
%%  the original line is commented out
 % of 2 or less will have a 1 dimensional volume of at most 200,000 and a
of 2 or less will have a 1-dimensional volume of at most 200,000 and a
0-dimensional boundary volume of 200,000.
%%  TLS edit
%%  The next line is not a sentence.
Since both of these numbers are less than $\infty$, the current is an integral current.

An integral 1-current can have an infinite
number of circles as long as the sum of the circumferences of the
circles is finite. For example, the following is an integral
current. Begin by enumerating all the points in the unit square having
rational coordinates, $p_1 = (x_1, y_1), p_2 = (x_1,y_2)$, ... .  Now
let $C(x,r)$ be the circle centered at $x$ with radius $r$, oriented
counterclockwise. Finally, define \[T = \bigcup_{i=1}^\infty
C\left(p_i, \frac{1}{2^i}\right).\] We then have that $M(T) = 2\pi$
and $M(\partial T) = 0$, so $T$ is a ``nice'' integral 1-current.

In the case of $k$-currents in $\R^n$, $n > 2$, things can be even
more interesting since we can now have an infinite number of pieces
having boundary (``non-closed'' pieces) as long as the total boundary
mass is finite.

\subsection*{Acknowledgment}

I am grateful to Thomas Scofield who caught many typos and made
suggestions that improved the exposition.

\bibliographystyle{plain}
\bibliography{/home/vixie/projects/templates/the_bib}

\end{document}